\documentclass[8pt]{amsproc}
\oddsidemargin=-1cm\evensidemargin=-1cm
\begin{document}
\numberwithin{equation}{section}

\def\1#1{\overline{#1}}
\def\2#1{\widetilde{#1}}
\def\3#1{\widehat{#1}}
\def\4#1{\mathbb{#1}}
\def\5#1{\frak{#1}}
\def\6#1{{\mathcal{#1}}}

\def\C{{\4C}}
\def\R{{\4R}}
\def\N{{\4N}}
\def\Z{{\4Z}}

\title{On a family of analytic discs attached to a real submanifold
$M\subset\mathbb{C}^{N+1}$}
\author{Valentin Burcea}
\begin{abstract} We construct a family of analytic discs attached  to a real submanifold
$M\subset\mathbb{C}^{N+1}$ of codimension $2$ near a CR
singularity. These discs are mutually disjoint and  form a smooth
hypersurface  $\widetilde{M}$  with boundary $M$ in a neighborhood
of the CR singularity.  As an application we prove that if $p$ is
a flat-elliptic CR singularity and if $M$ is nowhere minimal at
its CR points and does not contain a complex manifold of dimension
$(n-2)$, then $\widetilde{M}$ is a smooth Levi-flat hypersurface.
Moreover, if $M$ is real analytic we obtain that
$\widetilde{M}$ is real-analytic across the boundary manifold $M$.
\end{abstract}
\address{V. Burcea: School of Mathematics, Trinity College Dublin, Dublin 2, Ireland}
\email{valentin@maths.tcd.ie}
\thanks{This
project was supported by Science Foundation Ireland grant
 $10/RFP/MTH2878$.} \maketitle 


\def\Label#1{\label{#1}{\bf (#1)}~}


\def\cn{{\C^n}}
\def\cnn{{\C^{n'}}}
\def\ocn{\2{\C^n}}
\def\ocnn{\2{\C^{n'}}}


\def\dist{{\rm dist}}
\def\const{{\rm const}}
\def\rk{{\rm rank\,}}
\def\id{{\sf id}}
\def\tr{{\bf tr\,}}
\def\aut{{\sf aut}}
\def\Aut{{\sf Aut}}
\def\CR{{\rm CR}}
\def\GL{{\sf GL}}
\def\Re{{\sf Re}\,}
\def\Im{{\sf Im}\,}
\def\span{\text{\rm span}}
\def\Diff{{\sf Diff}}

\def\codim{{\rm codim}}
\def\crd{\dim_{{\rm CR}}}
\def\crc{{\rm codim_{CR}}}

\def\phi{\varphi}
\def\eps{\varepsilon}
\def\d{\partial}
\def\a{\alpha}
\def\b{\beta}
\def\g{\gamma}
\def\G{\Gamma}
\def\D{\Delta}
\def\Om{\Omega}
\def\k{\kappa}
\def\l{\lambda}
\def\L{\Lambda}
\def\z{{\bar z}}
\def\w{{\bar w}}
\def\Z{{\1Z}}
\def\t{\tau}
\def\th{\theta}

\emergencystretch15pt \frenchspacing

\newtheorem{Thm}{Theorem}[section]
\newtheorem{Cor}[Thm]{Corollary}
\newtheorem{Pro}[Thm]{Proposition}
\newtheorem{Lem}[Thm]{Lemma}

\theoremstyle{definition}\newtheorem{Def}[Thm]{Definition}

\theoremstyle{remark}
\newtheorem{Rem}[Thm]{Remark}
\newtheorem{Exa}[Thm]{Example}
\newtheorem{Exs}[Thm]{Examples}

\def\bl{\begin{Lem}}
\def\el{\end{Lem}}
\def\bp{\begin{Pro}}
\def\ep{\end{Pro}}
\def\bt{\begin{Thm}}
\def\et{\end{Thm}}
\def\bc{\begin{Cor}}
\def\ec{\end{Cor}}
\def\bd{\begin{Def}}
\def\ed{\end{Def}}
\def\br{\begin{Rem}}
\def\er{\end{Rem}}
\def\be{\begin{Exa}}
\def\ee{\end{Exa}}
\def\bpf{\begin{proof}}
\def\epf{\end{proof}}
\def\ben{\begin{enumerate}}
\def\een{\end{enumerate}}
\def\beq{\begin{equation}}
\def\eeq{\end{equation}}

\section{Introduction and Main Results}

Let $M\subset\mathbb{C}^{N+1}$ be a real submanifold.
A point $p\in M$ is called a CR singularity if
it is a  discontinuity point for the map $M\ni q\mapsto
\dim_{\mathbb{C}}T^{0,1}_{q}M$ defined near $p$. Here
$T_{q}^{0,1}M$ is the CR tangent space to $M$ at $q$.

We assume that $\codim_{\mathbb{R}}M=2$. Bishop considered the case when there
exist  $\left(z,w\right)$ holomorphical coordinates in $\mathbb{C}^{2}$,  such that near the CR
singularity $p=0$ the submanifold  $M\subset\mathbb{C}^{2}$ is
defined  by
 \begin{equation}
  w=z\overline{z}+\lambda\left(z^{2}+\overline{z}^{2}\right)+\mbox{O}(3),
  \label{1}\end{equation}
 where $\lambda\in\left[0,\infty\right)$ is a holomorphic invariant
 called the Bishop invariant. In the case when
 $\lambda\in\left[0,\frac{1}{2}\right)$, Kenig-Webster
 proved in ~\cite{KW1} the existence of a  unique family of $1$-dimensional analytic disks  shrinking to the CR singularity $p=0$. The real-analytic case was studied by Huang-Krantz in ~\cite{HY5}.

 Let $\left(z_{1},\dots,z_{N},w\right)$ be the coordinates from $\mathbb{C}^{N+1}$. In this paper, we consider the higher dimensional  case of $(\ref{1})$  when the submanifold $M\subset\mathbb{C}^{N+1}$ is defined near $p=0$ by
\begin{equation}
w=z_{1}\overline{z}_{1}+\lambda\left(z^{2}_{1}+\overline{z}^{2}_{1}\right)+Q\left(z_{1},\overline{z}_{1},z_{2},\overline{z}_{2},\dots,z_{N},\overline{z}_{N}\right)
+\mbox{O}(3),\label{1111}
\end{equation}
 where   $\lambda\geq 0$ and $Q\left(z_{1},\overline{z}_{1},z_{2},\overline{z}_{2},\dots,z_{N},\overline{z}_{N}\right)$ is a quadratic form not containing multiples of $z_{1}\overline{z}_{1}$, $z^{2}_{1}$ and $\overline{z}^{2}_{1}$.

In this paper, we extend  Kenig-Webster's Theorem from \cite{KW1}. We prove the following result:

\bt Let $M\subset\mathbb{C}^{N+1}$ be a smooth
 submanifold  defined locally near $p=0$ by (\ref{1111}) such that $\lambda$ is
 elliptic. Then there exists a family of
 regularly embedded analytic discs with boundaries on $M$ that are mutually disjoint and that forms a
smooth hypersurface $\widetilde{M}$ with boundary $M$ in a
neighborhood of the CR singularity $p=0$. \et

Then  $\widetilde{M}$ given  by Theorem $1.1$ is not necessary a Levi-flat hypersurface as in
Kenig-Webster's case from \cite{KW1} in $\mathbb{C}^{2}$.

The  existence problem of a Levi-flat hypersurface with prescribed boundary $S$ in $\mathbb{C}^{N+1}$ with $N\geq2$, was studied by Dolbeault-Tomassini-Zaitsev in ~\cite{DTZ}  under the following natural assumptions
on $S$:

(i) $S$ is compact, connected and  nowhere minimal at its CR
points;

(ii) $S$ does not contain a complex submanifold of dimension
$(n-2)$;

(iii) $S$ contains a finite number of flat elliptic CR
singularities.

We would like to mention that properties of  nowhere minimal CR submanifolds  were studied by Lebl \cite{L1}.

The  CR singularity $p=0$ is called elliptic if the quadratic part
from (\ref{1111}) is positive definite. We say that $p=0$ is a
''flat'' if  Definition $2.1$ from ~\cite{DTZ} is
satisfied. Under the precedent natural assumptions,  Dolbeault-Tomassini-Zaitsev proved the existence of a (possibly singular) Levi-flat hypersurface which
bounds $S$ in the sense of currents (see Theorem $1.3$,
~\cite{DTZ}).

The graph case was studied by  Dolbeault-Tomassini-Zaitsev in ~\cite{DTZ1}: Let
$\mathbb{C}^{N+1}=\left(\mathbb{C}^{N}_{z}\times\mathbb{R}_{u}\right)\times
\mathbb{R}_{v}$, where  $w=u+iv$, and let $\Omega$ be a bounded
strongly convex domain of $\mathbb{C}^{N}_{z}\times\mathbb{R}_{u}$
with smooth boundary $b\Omega$. Let $S\subset\mathbb{C}^{N+1}$,
$n\geq3$, be the graph of a function $g:b\Omega\longrightarrow
\mathbb{R}_{v}$ such that $S$ satisfies the natural assumptions (i), (ii), (iii). Under these assumptions, Dolbeault-Tomassini-Zaitsev proved the following result

 \bt Let $q_{1},q_{2}\in b\Omega$ be the projections of the complex
 points $p_{1},p_{2}$ of $S$, respectively. Then, there exists a
 Lipschitz function $f:\overline{\Omega}\longrightarrow\mathbb{R}_{v}$ which is smooth on
 $\overline{\Omega}-\left\{q_{1},q_{2}\right\}$ and such that
 $f|_{b\Omega}=g$ and $N=\mbox{graph}\left(f\right)-S$ is a
 Levi-flat hypersurface of $\mathbb{C}^{N+1}$. Moreover, each
 complex leaf of $M_{0}$  is the graph of a holomorphic function
 $\phi:\Omega'\longrightarrow\mathbb{C}$ where
 $\Omega'\subset\mathbb{C}^{n-1}$ is a domain with smooth boundary
 (that depends on the leaf) and $\phi$ is smooth on $\Omega'$.
  \et
As an application of Theorem $1.1$, we solve an open problem regarding the regularity of $f$ given by Theorem $1.2$ at $q_{1},q_{2}$,  proposed by  Dolbeault-Tomassini-Zaitsev in ~\cite{DTZ1}. By combining  Theorem $1.1$ and Theorem $1.2$,  we obtain the following result

\bt Let $M\subset\mathbb{C}^{N+1}$ be a smooth
 submanifold as in Theorem 1.2. Suppose $p$ is a point in $M$ such that $M$ is   defined  near $p=0$ by (\ref{1111}) satisfying the conditions that (i) $p=0$ is a flat-elliptic CR singularity (ii) any CR point of $M$ near $p=0$ is non-minimal, and (iii) $M$ does not contain a complex submanifold of  dimension $n-2$. Then $\widetilde{M}$ constructed by Theorem 1.1 is a
smooth Levi-flat hypersurface  with boundary $M$ in a neighborhood of $p=0$.
\et

In the real analytic case, our smoothness result combined with a similar  argument as in the papers \cite{HYA},\cite{HYB}  of Huang-Yin concerning the analyticity of the local hull of holomorphy,  gives the following result:

\bt Let $M\subset\mathbb{C}^{N+1}$ be a real analytic submanifold defined near $p=0$ by (\ref{1111}) and  that satisfies the assumptions  of Theorem $1.3$.
Then  $\widetilde{M}$ is a
Levi-flat hypersurface  real-analytic across the boundary manifold
$M$. \et

 We prove our results by following the lines of   Huang~\cite{HY4}, Kenig-Webster~\cite{KW1}, ~\cite{KW2} and in particulary the construction of analytic discs developed  by Huang-Krantz ~\cite{HY5}.  First, we  make
a perturbation along the CR singularity and then we find a holomorphic change of coordinates depending smoothly on a parameter. Then, we will adapt the methods used in
$\mathbb{C}^{2}$ by Huang-Krantz and Kenig-Webster in our case.  We would like to mention that versions of our result were  obtained in a higher codimensional case by Huang ~\cite{HY4} and Kenig-Webster~\cite{KW2}.

  The study of  CR singular real submanifolds in the complex space requires different methods than the case of CR manifolds.  We would like to mention here as a related topic the problem of finding a  normal form for a real submanifold in $\mathbb{C}^{2}$ defined by $(1.1)$. For $\lambda$ non-exceptional, this problem was solved by Moser-Webster in \cite{MW}. When $\lambda>\frac{1}{2}$, Gong constructed in \cite{G1}  an example of a real-analytic surface that is formally equivalent, but not biholomorphically equivalent to a quadric in $\mathbb{C}^{2}$. The case $\lambda=0$  was initiated by Moser in \cite{Mos} and completely understood  by Huang-Yin in the nice paper \cite{HY1}.  We would like to mention here also the paper \cite{HY2} of Huang-Yin where Moser's Theorem from \cite{Mos} is generalized  to the higher dimensional case.

$\bf{Acknowledgements.}$ I am grateful to my supervisor Prof. Dmitri Zaitsev for useful conversations and for bringing   this problem
to my attention. I would like thank  the  referee for helpful comments on the previous version of the manuscript.

\section{Preliminaries}
 \subsection{A Perturbation Along the CR singularity} Let $\Delta$ be the unit open disc from
$\mathbb{C}$ and let $\mbox{S}^{1}$ be its boundary. A map
$f:\overline{\Delta}\longrightarrow\mathbb{C}^{N+1}$ is called an
analytic disc if $f|_{\overline{\Delta}}$ is continuous and
$f|_{\Delta}$ is nonconstant and holomorphic. We say that  $f$ is an analytic disc
attached to $M$ if $f\left(\mbox{S}^{1}\right)\subset M$. 

We construct analytic discs attached to $M$  depending smoothly on
\begin{equation}X=\left(z_{2},\dots, z_{N}\right)=\left(x_{2}+iy_{2},\dots,x_{N}+iy_{N}\right)
\approx 0\in\mathbb{C}^{N-2}.\end{equation}

 By using the
notation $z=z_{1}$, our manifold $M$ is defined near $p=0$ by
\begin{equation}
w=z\overline{z}+\lambda\left(z^{2}+\overline{z}^{2}\right)+Q\left(z_{1},
\overline{z}_{1},z_{2},\overline{z}_{2},\dots,z_{N},\overline{z}_{N}\right)+\mbox{O}(3),\label{11}
\end{equation} or equivalently by
\begin{equation}   w= H_{0,0}\left(X\right)+\overline{z}H_{0,1}\left(X\right)+z
H_{1,0}\left(X\right)+z\overline{z}\left(1+H_{1,1}\left(X\right)\right)+\left(\lambda
+H_{2,0}\left(X\right)\right)z^{2}+\left(\lambda+H_{0,2}\left(X\right)\right)\overline{z}^{2}
+\mbox{O}\left(\left|z\right|^{3}\right),
\label{p1}
\end{equation}
where $H_{0,0}\left(X\right)$, $H_{1,0}\left(X\right)$,
$H_{0,1}\left(X\right)$, $H_{1,1}\left(X\right)$,
$H_{2,0}\left(X\right)$, $H_{0,2}\left(X\right)$ are smooth
functions vanishing at $X=0$.

 We prove  the following lemma:
 \bl Let
$M\subset\mathbb{C}^{2}$ be a real smooth submanifold defined near
$p=0$ by  $w=az+b\overline{z}+\mbox{O}\left(|z|^{2}\right)$. Then
\begin{equation}T_{0}^{c}M\neq\emptyset\Longleftrightarrow b=0.\end{equation}
 \el
\begin{proof}We need to solve the equations $\d f=\overline{\d}f=0$
at the point $z=w=0$. We compute:
\begin{equation}\begin{split}\d f|_{0}=\frac{\d f}{\d z}(0)dz+\frac{\d
f}{\d w}\left(0\right)dw=-dw+a d z,\quad \overline{\d}
f|_{0}=\frac{\d f}{\d \overline{z}}(0)d\overline{z}+\frac{\d f}{\d
\overline{w}}(0)d\overline{w}=b
d\overline{z}.\end{split}\end{equation} We obtain $adz=dw$ and $b
d\overline{z}=0$. It follows that $p=0$ is a CR singularity if and
only if $b=0$.
\end{proof}

We  make  a change of coordinates depending smoothly on $X\approx
0\in\mathbb{C}^{N-2}$  preserving  the CR singularity $p=0$ :

\bp There exists a biholomorphic change of coordinates in
$\left(z,w\right)$ depending smoothly on
$X\approx0\in\mathbb{C}^{N-2}$ that sends (\ref{p1}) to a
submanifold defined by
\begin{equation}w=z\overline{z}+\lambda\left(X\right)\left(z^{2}+\overline{z}^{2}\right)+\mbox{O}\left(|z|^{3}\right),\label{p8}
\end{equation}
preserving the CR singularity $p=0$. Here
$0\leq\lambda\left(X\right)<\frac{1}{2}$ for $X\approx
0\in\mathbb{C}^{N-2}$ and $\lambda\left(0\right)=\lambda$.
\ep
\begin{proof} We consider a local defining function for $M$ near $p=0$
\begin{equation} 
f\left(z,X,w\right)=-w+H_{0,0}\left(X\right)+\overline{z}H_{0,1}\left(X\right)+z
H_{1,0}\left(X\right)+z\overline{z}\left(1+H_{1,1}\left(X\right)\right) +\left(\lambda
+H_{2,0}\left(X\right)\right)z^{2}+\left(\lambda+H_{0,2}\left(X\right)\right)\overline{z}^{2}
+\mbox{O}\left(\left|z\right|^{3}\right).\label{df}
 \end{equation}
 
Each fixed $X\approx 0\in\mathbb{C}^{N-2}$ defines us a real
submanifold in $\mathbb{C}^{2}$ which may not have a CR
singularity at the point $z=w=0$ because $H_{0,1}\left(X\right)$
may be different than $0$ (see Lemma $2.1$). Therefore we need to
make a change of coordinates in $\left(z,w\right)$ depending
smoothly on $X\approx0\in\mathbb{C}^{N-2}$ that perturbs the CR
singularity $p=0$. We consider the following equation
\begin{equation}0=\frac{\d f}{\d\overline{z}}=H_{0,1}\left(X\right)+\left(1+H_{1,1}\left(X\right)\right)z
+B\left(z,\overline{z},X\right),\label{9000}\end{equation} where
$B\left(z,\overline{z},X\right)$ is a smooth function. 

 Because
$H_{1,1}\left(0\right)=0$, by applying the implicit function
theorem,  we obtain a smooth solution $z_{0}=z_{0}\left(X\right)$
for (\ref{9000}). By making  the translation
$\left(w',z'\right)=\left(w,z+z_{0}\left(X\right)\right)$, the
equation (\ref{p1}) becomes
\begin{equation} w=z
C_{1,0}\left(X\right)+z\overline{z}\left(1+C_{1,1}\left(X\right)\right)+\left(\lambda
+C_{2,0}\left(X\right)\right)z^{2}+(\lambda+C_{0,2}\left(X\right))\overline{z}^{2}+\mbox{O}\left(\left|z\right|^{3}\right),\label{90}\end{equation}
where $C_{1,0}\left(X\right)$, $C_{1,1}\left(X\right)$,
$C_{2,0}\left(X\right)$, $C_{0,2}\left(X\right)$ are smooth
functions vanishing at $X=0$. Let
$\gamma\left(X\right)=1+C_{1,1}\left(X\right)$,
$\Lambda_{1}\left(X\right)=\lambda+C_{2,0}\left(X\right)$,
$\Lambda_{2}\left(X\right)=\lambda+C_{0,2}\left(X\right)$. In the
new coordinates
$\left(w,z\right):=\left(\left(w-C_{1,0}\left(X\right)z\right)/\gamma\left(X\right),z\right)$,
the equation (\ref{90}) becomes
\begin{equation}w=z\overline{z}+\Lambda_{1}\left(X\right)z^{2}+\Lambda_{2}\left(X\right)\overline{z}^{2}
+\mbox{O}\left(\left|z\right|^{3}\right).\end{equation} Next, we
consider a map $\Theta\left(X\right)$ such that
$\Lambda_{2}\left(X\right)e^{-2i\Theta\left(X\right)}\geq 0$.
Changing the coordinates
 $\left(w,z\right):=\left(w,ze^{i\Theta\left(X\right)}\right)$, we can assume $\Lambda_{2}\left(X\right)\geq 0$.
 Changing again the coordinates
$\left(w,z\right):=\left(w+\left(\Lambda_{1}\left(X\right)-\Lambda_{2}\left(X\right)\right)z^{2},z\right)$
we obtain (\ref{p8}).
\end{proof}

We  write
\begin{equation}M:
w=z\overline{z}+\lambda\left(X\right)\left(z^{2}+\overline{z}^{2}\right)+P\left(z,X\right)+iK\left(z,X\right),\label{ecM}\end{equation}
where $P\left(z,X\right)$ and $K\left(z,X\right)$ are real smooth
functions. We prove an extension of Lemma $1.1$ from ~\cite{KW1}:

\bp There exists a holomorphic change of coordinates in
$\left(z,w\right)$ depending smoothly on
$X\approx0\in\mathbb{C}^{N-2}$ in which $K$
 and its partial derivatives in $z$ and $\overline{z}$ of order less or equal to $l$ vanish at $z=0$.
 \ep
\begin{proof}By making the substitution
$\left(z'\left(X\right),w'\left(X\right)\right)=\left(z,w+B\left(z,X,w\right)\right)$
and by (\ref{ecM}) it follows that
\begin{equation}M:w'=q\left(z,X\right)+P\left(z,X\right)+i K\left(z,X\right)+\Re
B\left(z,X,w\right)+i\Im
B\left(z,X,w\right),\label{212}\end{equation} where
$q\left(z,X\right)=z\overline{z}+\lambda\left(X\right)\left(z^{2}+\overline{z}^{2}\right)$.

We want to make the derivatives in $z$ of order less than $l$ of
$i\left(K\left(z,X\right)+\Im B\left(z,X,w\right)\right)$ vanish
at $z=0$. By multiplying (\ref{212}) by $i=\sqrt{-1}$, our problem is
reduced to the following general equation
\begin{equation}\Re B\left(z,X,q\left(z,X\right)+P\left(z,X\right)+iK\left(z,X\right)\right)
=f\left(z,\overline{z},X\right),\label{40}\end{equation} where
$f\left(z,\overline{z},X\right)$ is a real formal power series in
$\left(z,\overline{z},X\right)$ with cubic terms in $z$ and
$\overline{z}$ with coefficients depending smoothly on
$X\approx0\in\mathbb{C}^{N-2}$.

 We write as follows
\begin{equation}\begin{split}&f\left(z,\overline{z},X\right)
=\displaystyle\sum_{m=3}^{l}f_{m}\left(z,\overline{z},X\right),
\hspace{0.1
cm}f_{m}\left(z,\overline{z},X\right)=\displaystyle\sum_{j_{1}+j_{2}=m}c_{j_{1},j_{2}}^{m}\left(X\right)
z^{j_{1}}\overline{z}^{j_{2}},\hspace{0.1
cm}c_{j_{1},j_{2}}^{m}\left(X\right)=\overline{c_{j_{2},j_{1}}^{m}\left(X\right)},\\&\quad\quad\quad
B\left(z,X,w\right)=\displaystyle\sum_{m=3}^{l}B_{m}\left(z,X,w\right),\quad
B_{m}\left(z,X,w\right)=\displaystyle\sum_{j_{1}+2
j_{2}=m}b_{j_{1},j_{2}}^{m}\left(X\right)z^{j_{1}}w^{j_{2}}.
\end{split}
\end{equation}

We  solve inductively (\ref{40}) by using the following lemma:
 \bl The equation (\ref{40}) has a 
unique solution with the normalization condition $\Im
B_{m}\left(0,X,u\right)=0$.
 \el
\begin{proof}
We define the weight of $z$ to be $1$ and the weight of $w$ to be
$2$. We say that the polynomial $B_{m}\left(z,X,w\right)$ has
weight $m$ if
$B_{m}\left(tz,X,t^{2}w\right)=t^{m}B_{m}\left(z,X,w\right)$. Let
$\mathbb{B}_{m}$ be the space of all such homogeneous holomorphic
polynomials in $\left(z,w\right)$ of weight $m$ satisfying the
normalization condition with coefficients depending smoothly on
$X\approx0\in\mathbb{C}^{N-2}$ and let $\mathbb{F}_{m}$ be the
space of all homogeneous polynomials
$f_{m}\left(z,\overline{z},X\right)$ of bidegree
$\left(k,l\right)$ in $\left(z,\overline{z}\right)$ with $k+2l=m$
with coefficients depending smoothly on
$X\approx0\in\mathbb{C}^{N-2}$.

 We can rewrite (\ref{40}) as
follows
\begin{equation}B_{m}\left(z,X,q\left(z,X\right)+P\left(z,X\right)
+iK\left(z,X\right)\right)=B_{m}\left(z,X,q\left(z,X\right)\right)
+\mbox{O}\left(\left|z\right|^{m+1}\right)
.\label{400}\end{equation}

 In order to solve (\ref{400}), it is
enough to prove that we have a linear invertible transformation
\begin{equation}\varphi\left(X\right):\mathbb{B}_{m}\ni
B_{m}\left(z,X,w\right)\mapsto \Re
B_{m}\left(z,X,q\left(z,X\right)\right)\in
\mathbb{F}_{m},\end{equation} depending smoothly on $X\approx
0\in\mathbb{C}^{N-2}$. By Lemma $1.1$ from the paper ~\cite{KW1} of Kenig-Webster, it follows
that $\varphi\left(X\right)$ is invertible for
$X=0\in\mathbb{C}^{N-2}$. By   continuity, it follows that
$\varphi\left(X\right)$ is invertible. If it is necessary, then we
shrink the range of $X\approx 0\in\mathbb{C}^{N-2}$.
\end{proof}

The proof is completed now by induction and by using Lemma $2.4$.
\end{proof}

\subsection{Preliminary Preparations} Let $w=u+iv$ and $I_{\epsilon}:=\left(-\epsilon,\epsilon\right)\subset\mathbb{R}$, for $0<\epsilon<<1$. We assume that $M$ is defined by (\ref{ecM})
and satisfies the properties of Proposition $2.3$.

In order to define  a family of attached discs to the manifold
$M$, we define the following domain
\begin{equation}\begin{split}D_{X,r}=\left\{ z\in\mathbb{C};\hspace{0.1 cm}v=0,\hspace{0.1 cm}
q\left(z,X\right)+P\left(z,X\right)\leq
u<\epsilon\right\},\label{p4}
\end{split}\end{equation}
where $u=r^{2}$. By similar arguments as in the paper ~\cite{HY4} of Huang, it follows that $D_{X,r}$  is a simply connected bounded set of $\mathbb{C}$. Therefore
there exists a  unique mapping $r\sigma_{X,r}:\Delta\rightarrow
D_{X,r}$ such that $\sigma_{X,r}\left(0\right)=0$ and
$\sigma'_{X,r}\left(0\right)>0$. Then, for $0<r<<1$  we can define
the following family of curves depending smoothly on
$X\approx0\in\mathbb{C}^{N-2}$
\begin{equation}\gamma_{X,r}=\left\{z\in\mathbb{C};\hspace{0.1
cm}q\left(z,X\right)+P\left(z,X\right)=r^{2}\right\}.
\label{gmma}\end{equation} Next,  we define the following family
of analytic discs
\begin{equation}\left\{\left(r\sigma_{X,r},X,r^{2}\right)\right\}_{X\approx0\in\mathbb{C}^{N-2},\hspace{0.1 cm}0<r<<1}.\label{p2}
\end{equation}
The family of analytic discs shrinks to
$\left\{0\right\}\times\mathcal{O}\times\left\{0\right\}$ as
$r\mapsto0$, where $0\in\mathcal{O}\subset\mathbb{C}^{N-2}$ and
fills up the following domain
\begin{equation}\widetilde{M}_{0}=\left\{\left(z,X,u\right)\in\mathbb{C}
\times\mathbb{C}^{N-2}\times\mathbb{R};\hspace{0.1cm}\left\|X\right\|<<1,
\hspace{0.1 cm}q\left(z,X\right)+P\left(z,X\right)\leq
u\right\}\label{m0} .\end{equation}

\subsection{The Hilbert Transform on a Variable Curve} Let
$\gamma_{X,r}$ given by (\ref{gmma}), where $r$ is taken very
small. For a function $\varphi_{X,r}\left(\theta\right)$ defined
on $\gamma_{X,r}$ we define its Hilbert transform
$\mbox{H}_{X,r}\left [\varphi_{X,r}\right]$ to be the boundary
value of a function holomorphic inside $\gamma_{X,r}$, with its
imaginary part vanishing at the origin. For more informations about Hilbert's transform we mention here the book \cite{Hb} of Helmes.

For $\alpha\in \left(0,1\right)$, we define the following Banach
spaces:
\begin{equation}\begin{split}&\mathcal{C}^{\alpha}:=\left\{ u:\gamma_{X,r}\longrightarrow
\mathbb{R};\hspace{0.1 cm}\left\|u\right\|_{\alpha}:=\displaystyle\sup_{x
\in \gamma_{X,r}}\left|u(x)\right| +\displaystyle\sup _{ x,y \in
\gamma_{X,r} \atop x\neq y}\frac{\left|u(x)-u(y)\right|}{\left|x-y\right|^{\alpha}}<\infty
\right\}, \\& 
\mathcal{C}^{k,\alpha}:=\left\{ u:\gamma_{X,r}\longrightarrow
\mathbb{R};\hspace{0.1
cm}\left\|u\right\|_{k,\alpha}:=\displaystyle\sum
_{\left|\beta\right|\leq
k}\left\|\mbox{D}^{\beta}u\right\|_{\alpha}<\infty \right\}.
\end{split}\end{equation}

 Let  $X:=\left\{x_{2},y_{2},\dots,x_{N},y_{N}\right\}$. The following result can be proved by using the
same lines as in Kenig-Webster's paper ~\cite{KW1}  (Theorem $2.5$) or from Kenig-Webster's paper ~\cite{KW2}:
 \bp As $r\rightarrow 0$ and $X\approx 0\in\mathbb{C}^{N-2}$ we have
 \begin{equation}\begin{split}\left\|\mathcal{H}_{X,r}\right\|_{j,\alpha}=\mbox{O}\left(r\right),
 \hspace{0.1
 cm}\mbox{for all $j\leq l-2$};\quad\left\|\left(\d_{X}^{\left|
I\right|}\d^{s}_{r}\right)\mathcal{H}_{X,r}\right\|_{j,\alpha}=\mbox{O}(1)
 ,\hspace{0.1
cm}\mbox{for all $j+2s\leq l-4$},\hspace{0.1
cm}I\in\mathbb{N}^{N-2}.\end{split}\end{equation}
 \ep

\subsection{An Implicit Functional Equation} During this section we
work in the Holder space
$\left(\mathcal{C}^{j,\alpha},\left\|\cdot\right\|_{j,\alpha}\right)$.
We employ ideas developed by Huang-Krantz in ~\cite{HY5}, Huang in
 ~\cite{HY4},  Kenig-Webster in ~\cite{KW1},~\cite{KW2} and we
define the following auxiliary hypersurface
\begin{equation}M_{0}=\left\{\left(z,X,u\right)\in\mathbb{C}\times\mathbb{C}^{N-2}\times\mathbb{R};\hspace{0.1 cm}
\left\|X\right\|<<1,\hspace{0.1
cm}q\left(z,X\right)+P\left(z,X\right)=u<\epsilon\right\},\label{m00}\end{equation}
where $\epsilon>0$ is small enough and $w=u+iv$. We would like to
find a map of the following type
\begin{equation}T=T\left[X\right]:=\left(z\left(1+\mathcal{F}\left(z,X,r\right)\right)
,\mathcal{B}\left(z,X,r\right)\right)\label{p3}\end{equation} such
that $T\left(M_{0}\right)\subseteq M$. Here $\mathcal{F}$,
$\mathcal{B}$ are holomorphic functions in $z$ and smooth in
$\left(X,r\right)$. It follows that
\begin{equation}\mathcal{B}\left(z,X,r\right)|_{\gamma_{X,r}}=\left(q+P+iK\right)
\left(z+z\mathcal{F}\left(z,X,r\right),X\right)|_{\gamma_{X,r}},\label{224}\end{equation}
where $\gamma_{X,r}$ is the curve defined by (\ref{gmma}). By
using the Hilbert transform on the curve $\gamma_{X,r}$ and by
dividing by $r^{2}$ the equation (\ref{224}), it follows that
there exists a smooth function $V\left(X,r\right)$ such that
\begin{equation} q\left(z\left(1+\mathcal{F}\left(z,X,r\right)\right),X\right)|_{\gamma_{X,r}}=
-P\left(z\left(1+\mathcal{F}\left(z,X,r\right)\right),X\right)|_{\gamma_{X,r}}
 -\mathcal{H}_{X,r}\left[K\left(z\left(1+\mathcal{F}\left(z,X,r\right)\right)
,X\right)\right]|_{\gamma_{X,r}}+V\left(X,r\right).\label{ecfu} \end{equation}

We follow Huang-Krantz's strategy from ~\cite{HY5} and we  define
 the following functional
\begin{equation} \Omega\left(\mathcal{F},X,r\right)=\frac{q\left(z\left(1+\mathcal{F}\right),X\right)
+P\left(z\left(1+\mathcal{F}\right),X\right)}{r^{2}}|_{\gamma_{X,r}}
,\label{1004} \end{equation} where
$\mathcal{F}=\mathcal{F}\left(z,X,r\right)$.

 By linearizing in
$\mathcal{F}=0$ the functional defined in (\ref{1004}), the
equation (\ref{ecfu})
 becomes
\begin{equation}1+\Omega'\left(\mathcal{F},X,r\right)+\Omega_{1}\left(\mathcal{F},X,r\right)
+\frac{1}{r^{2}}\mathcal{H}_{X,r}\left[K\left(z\left(1+\mathcal{F}\right),X\right)\right]|_{\gamma_{X,r}}-\frac{V\left(X,r\right)}{r^{2}}=0,\label{ecuatiefuntionala}\end{equation}
where $\mathcal{F}=\mathcal{F}\left(z,X,r\right)$ and
$\Omega_{1}\left(\mathcal{F}\left(z,X,r\right),X,r\right)$, are
terms that are coming from the Taylor expansion of
 $P\left(z,X\right)$ and
\begin{equation}\Omega'\left(\mathcal{F},X,r\right)=\frac{2}{r^{2}}
\Re\left\{\left(q+P\right)_{z}\left(z,X\right)z\mathcal{F}\right\}|_{\gamma_{X,r}}.\end{equation}

We put the normalization condition $V\left(X,r\right)=r^{2}$. In
order to find a solution $\mathcal{F}$ in the Holder space
$\left(\mathcal{C}^{j,\alpha},\|\cdot\|_{j,\alpha}\right)$ for
(\ref{ecuatiefuntionala}), we need to study the regularity
properties of the functional $\Omega$. We consider the following
notation
\begin{equation}\begin{split}\mathcal{C}_{X,r}\left(z\right)=\frac{2}{r^{2}}\Re\left\{
\left(q+P\right)_{z}\left(z,X\right)z\right\}|_{\gamma_{X,r}}.
\end{split}\end{equation}

Because $\mathcal{C}_{X,r}\left(z\right)\neq 0$ for $\left|r \right|<<1$,
$X\approx 0\in\mathbb{C}^{N-2}$, we can write
 $\mathcal{C}_{X,r}\left(z\right)=\mathcal{A}\left(z,X,r\right)\mathcal{B}\left(z,X,r\right)$ with
\begin{equation}\mathcal{A}\left(z,X,r\right)=\left|\mathcal{C}_{X,r}\left(z\right)\right|,\quad
\mathcal{B}\left(z,X,r\right)=\frac{\mathcal{C}_{X,r}\left(z\right)}{\left|\mathcal{C}_{X,r}\left(z\right)\right|}.\end{equation}

Then $\ln \mathcal{B}\left(z,X,r\right)$ is a well-defined smooth
function in $\left(z, X,r\right)$. Among the lines developed by Huang-Krantz in ~\cite{HY5},
we define the following function
\begin{equation}\mathcal{C}^{\star}\left(z,X,r\right)=\frac{e^{i\mathcal{H}_{X,r}\left(\ln
\mathcal{B}\left(z,X,r\right)\right)}}{\mathcal{A}\left(z,X,r\right)}.\end{equation}
Then $\mathcal{C}^{\star}$ is a smooth positive function  and
$D\left(z,X,r\right):=\mathcal{C}^{\star}\left(z,X,r\right)\mathcal{C}\left(z,X,r\right)$
 is holomorphic  in $z$, smooth in $\left(X,r\right)$. We write  $D\left(z,X,r\right)\mathcal{F}\left(z,X,r\right)\equiv
U\left(z,X,r\right)+\sqrt{-1}\mathcal{H}_{X,r}\left[U\left(z,X,r\right)\right]$.
Because $D\left(z,X,r\right)\neq 0$,
 we can rewrite (\ref{ecuatiefuntionala}) as follows
\begin{equation}\begin{split}
U\left(z,X,r\right)=&-C^{\star}\left(z,X,r\right)\left(
\Omega_{1}\left(\frac{U\left(z,X,r\right)+i\mathcal{H}_{X,r}
\left[U\left(z,X,r\right)\right]}{D\left(z,X,r\right)},X,r\right)
\right) \\&-C^{\star}\left(z,X,r\right)
\frac{1}{r^{2}}\mathcal{H}_{X,r}\left[K\left(z\left(1+\frac{U\left(z,X,r\right)+i\mathcal{H}_{X,r}
\left[U\left(z,X,r\right)\right]}{D\left(z,X,r\right)}\right),X\right)\right].\end{split}
 \label{Uec}\end{equation}

 We summarize all the precedent  computations and we obtain the
 following regularity result
 \bt The equation (\ref{Uec}) has a  unique solution in the Banach space
$\left(\mathcal{C}^{j,\alpha},\hspace{0.1
cm}\|\cdot\|_{j,\alpha}\right)$ such that
\begin{equation}\begin{split}\left\|U\right\|_{j,\alpha}=\mbox{O}\left(r^{l-2}\right),
\hspace{0.1 cm}\mbox{for all}\hspace{0.1 cm} j\leq l-2;\quad
\left\|\left(\d_{X}^{\left|I\right|}\d_{r}^{s}\right)U\right\|_{j,\alpha}=\mbox{O}\left(r^{l-s-2}\right)
,\hspace{0.1 cm}\mbox{for all $j+2s\leq l-4$},\hspace{0.1
cm}I\in\mathbb{N}^{N-2}.
\end{split}
\end{equation}
\et
\begin{proof} The solution $U$ and its uniqueness follows by applying the
implicit function theorem. We denote by
$\Lambda_{1}\left(U,X,r\right)$ and
$\Lambda_{2}\left(U,X,r\right)$ the first  and the second term   from (\ref{Uec}). It follows that $$\left\|U\right\|_{j,\alpha}\leq
\left\|\Lambda_{1}\left(U,X,r\right)\right\|_{j,\alpha}+
 \left\|\Lambda_{2}\left(U,X,r\right)\right\|_{j,\alpha}\leq
 \left\|\Lambda_{1}\left(U,X,r\right)\right\|_{j,\alpha}+\mbox{O}\left(r^{l-2}\right)\leq C \left\|U\right\|^{2}_{j,\alpha}+\mbox{O}\left(r^{l-2}\right),$$ for some $C>0$. It follows
 that $\left\|U\right\|_{j,\alpha}=\mbox{O}\left(r^{l-2}\right)$.

 The proof of the second regularity property goes after the
 previous line. Differentiating  with $r$ the equation (\ref{Uec}) it follows
 that  $$\partial_{r}U=\partial_{r}\Lambda_{1}\left(U,X,r\right)
 +\partial_{U}\Lambda_{1}\left(U,X,r\right)\left[\partial_{r}U\right]
 +\partial_{r}\Lambda_{2}\left(U,X,r\right)
 +\partial_{U}\Lambda_{2}\left(U,X,r\right)\left[\partial_{r}U\right].$$
 
 By Proposition $2.3$ and Proposition $2.5$ we obtain that
 $\left\|\partial_{r}U\right\|_{j,\alpha}=\mbox{O}\left(r^{l-2-1}\right)$.
Because $P\left(z,X\right)=\mbox{O}\left(z^{3}\right)$ and
$K\left(z,X\right)=\mbox{O}\left(z^{l}\right)$, by taking higher derivatives of $x$ in (\ref{Uec}) it follows that the
differentiation of any order with $x\in X$ does
 not affect the estimates. Therefore the
 second estimates follow immediately.
\end{proof}

We write that
\begin{equation}\mathcal{F}_{X,r}\left[\varphi_{X,r}\right]=\frac{U\left(z,X,r\right)+i\mathcal{H}_{X,r}
\left[U\left(z,X,r\right)\right]}{D\left(z,X,r\right)}:=\varphi_{X,r}+i\mathcal{H}_{X,r}\left[\varphi_{X,r}\right]
,\label{1001}\end{equation} where $\left\|\varphi_{X,r}
\right\|_{j,\alpha}=\mbox{O}\left(r^{l-2}\right)$, for all $j\leq
l-2$ and
$\left\|\left(\d_{X}^{\left|I\right|}\d_{r}^{s}\right)\varphi_{X,r}\right\|_{j,\alpha}=\mbox{O}\left(r^{l-s-2}\right)
,\hspace{0.1 cm}\mbox{for all $j+2s\leq l-4$},\hspace{0.1
cm}I\in\mathbb{N}^{N-2}$.

\section{A Family of Analytic Discs and Proofs of Main Results}
\subsection{A Family of Analytic Discs}

 We construct a continuous mapping $T$ defined on $\widetilde{M}_{0}$ into
$\mathbb{C}^{2}$ that is holomorphic in $z$ for each fixed
$u=r^{2}$ and that maps slice by slice the hypersurface $M_{0}$
into $M$. Let $\varphi_{X,r}$ be the function defined  by
(\ref{1001}). Then
\begin{equation}\begin{split}\quad\quad\mathcal{F}_{X,r}\left[\varphi_{X,r}\right]
=\varphi_{X,r}+i\mathcal{H}_{X,r}\left[\varphi_{X,r}\right], \quad
\mathcal{B}_{X,r}\left[\varphi_{X,r}\right]=\left(q+P+iK\right)
\left(z+z\mathcal{F}_{X,r}\left[\varphi_{X,r}\right],X\right).
\end{split}
\label{fd1}
\end{equation}
We extend these functions to  $\widetilde{M}_{0}$ by the Cauchy
integral as follows
\begin{equation}\begin{split}&\mathcal{F}\left(\zeta,X,r\right)=\mathcal{C}\left(\mathcal{F}_{X,r}
\left[\varphi_{X,r}\right]\right)(\zeta)\equiv\frac{1}{2\pi
i}\displaystyle\int_{0}^{2\pi}\frac{\mathcal{F}_{X,r}\left[\varphi_{X,r}\right]
\left(\theta\right)z_{\theta}\left(\theta,X,r\right)}{z\left(\theta,X,r\right)-\zeta}d\theta,
\\& \mathcal{B}\left(\zeta,X,r\right)=\mathcal{C}\left(\mathcal{B}_{X,r}
\left[\varphi_{X,r}\right]\right)(\zeta)\equiv\frac{1}{2\pi i}
\displaystyle\int_{0}^{2\pi}\frac{\mathcal{B}_{X,r}\left[\varphi_{X,r}\right]
\left(\theta\right)z_{\theta}\left(\theta,X,r\right)}{z\left(\theta,X,r\right)-\zeta}d\theta,
\end{split}
\label{fd2}
\end{equation}
where $z=z\left(\theta,X,r\right)$ is a parameterization of the
curve $\gamma_{X,r}$  defined by  (\ref{gmma}).

We define $T$ by (\ref{fd1}). Then $T$ is continuous by
construction up to the boundary on each slice $\left(X,r\right)$=constant.
In order to obtain the regularity of $T$, we have to
bound the derivatives in $\left(z,X,u\right)$ of $\mathcal{F}$ and
$\mathcal{B}$. We state the following lemma:

\bl For all $j+2s\leq l-4$, $I\in\mathbb{N}^{N-2}$ as $r\mapsto
0$, we have
\begin{equation}\d_{\theta}^{j}\d_{X}^{\left|I\right|}\d_{r}^{s}\mathcal{F}\left(z,X,r\right)
=\mbox{O}\left(r^{l-s-2}\right),\quad\d_{z}^{j}\d_{X}^{\left|I\right|}\d_{r}^{s}
\mathcal{B}\left(z,X,r\right)=\mbox{O}\left(r^{l-s}\right).\end{equation}
 \el
The proof of the predent Lemma follows by the lines
 of Lemma $4.1$ proof from Kenig-Webster's paper ~\cite{KW1}.

\bt Let $M$ defined by (\ref{ecM}) with
$P\left(z,X\right)=\mbox{O}\left(z^{3}\right)$,
$K\left(z,X\right)=\mbox{O}\left(z^{l}\right)$, $l\geq 7$, $T$
extended by (\ref{fd2}). Then
$\widetilde{M}=T\left(\widetilde{M}_{0}\right)$ is a complex
manifold-with-boundary regularly foliated by discs embedded of
class $\mathcal{C}^{\frac{l-7}{3}}$.
 \et
\begin{proof} Since
$\d_{u}=\frac{1}{2r}\d_{r}$, it follows that
\begin{equation}\d_{z}^{j}\d_{X}^{\left|I\right|}\d_{u}^{s}\mathcal{F}_{X,r}\left(z,X,r\right)=\mbox{O}\left(r^{l-2s-j-2}\right),
\quad\d_{z}^{j}\d_{X}^{\left|I\right|}\d_{u}^{s}\mathcal{F}_{X,r}\left(z,X,r\right)=\mbox{O}\left(r^{l-2s}\right),\end{equation}
and these derivatives remain bounded for all $j+2s\leq l-4$,
$I\in\mathbb{N}^{N-2}$. It follows  that the jacobian matrix $DT$
of $T=T\left(X\right)$ is the identity matrix.
\end{proof}

\subsection{Proof of Theorem $1.1$} Let $M$, $\widetilde{M}$, $T$
as in Theorem $3.2$. Using the techniques from \cite{M},  \cite{NWY} together with  an extended reflection principle
as in  the paper ~\cite{KW2} of Kenig-Webster, we construct smooth extension of $T$ past every
point of $ M_{0}-\left\{0\right\}$. By similar arguments as in the papers
~\cite{KW1}, ~\cite{KW2} of Kenig-Webster, we obtain that $M\cup \widetilde{M}$ is a smooth
manifold-with-boundary $M$ in a neighborhood of the CR singular point $p=0$.

\subsection{Proof of Theorem $1.3$} Since the hypersurface given
by Theorem $2.1$ is Levi-flat, it follows each of our analytic
discs is a reparameterization of an analytic disc contained
inside. By dimension reasons, it follows that the under the
hypothesis of Theorem $1.2$, the hypersurfaces given by Theorem
$1.1$ and Theorem $1.2$ are the same.

\subsection{Proof of Theorem $1.4$} We can  study now the hull of $M$ near $p=0$ when $M$ is assumed
to be real-analytic. The hypersurface $M_{0}$ defined by
(\ref{m00}) is foliated by   the family of analytic discs defined
by (\ref{p2}) and therefore $\widetilde{M}$ is foliated by the family of
analytic discs defined by (\ref{fd2}). By similar arguments as in    the papers  \cite{HYA},\cite{HYB} of  Huang-Yin, we obtain our result. The author believes that the arguments from  the paper \cite{HY5} of Huang-Krantz  or from the paper \cite{HY4}  of Huang, can be adapted in order to  prove the analyticity in our case.

\end{document}